\documentclass[reqno,11pt]{article}
\usepackage{mathrsfs}
\usepackage{amssymb}
\usepackage[usenames,dvipsnames]{xcolor}
\usepackage{amsthm}
\usepackage{amsmath}
\usepackage{amsfonts}
\usepackage{enumerate}
\usepackage{txfonts}
\usepackage{bm}

\usepackage{graphicx}

\numberwithin{equation}{section}

\newtheorem{theorem}{Theorem}[section]
\newtheorem{lemma}[theorem]{Lemma}
\newtheorem{corollary}[theorem]{Corollary}
\newtheorem{proposition}[theorem]{Proposition}

\theoremstyle{definition}
\newtheorem{definition}[theorem]{Definition}
\newtheorem{assumption}[theorem]{Assumption}
\newtheorem{example}[theorem]{Example}

\theoremstyle{remark}
\newtheorem{remark}[theorem]{Remark}

\begin{document}

\def\be{\begin{eqnarray}}
\def\ee{\end{eqnarray}}
\def\p{\partial}
\def\no{\nonumber}
\def\e{\epsilon}
\def\de{\delta}
\def\De{\Delta}
\def\om{\omega}
\def\Om{\Omega}
\def\f{\frac}
\def\th{\theta}
\def\la{\lambda}
\def\lab{\label}
\def\b{\bigg}
\def\var{\varphi}
\def\na{\nabla}
\def\ka{\kappa}
\def\al{\alpha}
\def\La{\Lambda}
\def\ga{\gamma}
\def\Ga{\Gamma}
\def\ti{\tilde}
\def\wti{\widetilde}
\def\wh{\widehat}
\def\ol{\overline}
\def\ul{\underline}
\def\Th{\Theta}
\def\si{\sigma}
\def\Si{\Sigma}
\def\oo{\infty}
\def\q{\quad}
\def\z{\zeta}
\def\co{\coloneqq}
\def\eqq{\eqqcolon}
\def\di{\displaystyle}
\def\bt{\begin{theorem}}
\def\et{\end{theorem}}
\def\bc{\begin{corollary}}
\def\ec{\end{corollary}}
\def\bl{\begin{lemma}}
\def\el{\end{lemma}}
\def\bp{\begin{proposition}}
\def\ep{\end{proposition}}
\def\br{\begin{remark}}
\def\er{\end{remark}}
\def\bd{\begin{definition}}
\def\ed{\end{definition}}
\def\bpf{\begin{proof}}
\def\epf{\end{proof}}
\def\bex{\begin{example}}
\def\eex{\end{example}}
\def\bq{\begin{question}}
\def\eq{\end{question}}
\def\bas{\begin{assumption}}
\def\eas{\end{assumption}}
\def\ber{\begin{exercise}}
\def\eer{\end{exercise}}
\def\mb{\mathbb}
\def\mbR{\mb{R}}
\def\mbZ{\mb{Z}}
\def\mc{\mathcal}
\def\mcS{\mc{S}}
\def\ms{\mathscr}
\def\lan{\langle}
\def\ran{\rangle}
\def\lb{\llbracket}
\def\rb{\rrbracket}

\title{Decay properties of axially symmetric D-solutions to the steady Navier-Stokes equations}

\author{{Shangkun Weng\thanks{School of mathematics and statistics, Wuhan University, Wuhan, Hubei Province, China, 430072. Email: skweng@whu.edu.cn}}}

\pagestyle{myheadings} \markboth{Decay properties of the axially symmetric D-solutions to steady Navier-Stokes equations}{}\maketitle

\begin{abstract}
  We investigate the decay properties of smooth axially symmetric D-solutions to the steady Navier-Stokes equations. The achievements of this paper are two folds. One is improved decay rates of $u_{\th}$ and $\na {\bf u}$, especially we show that $|u_{\th}(r,z)|\leq c\left(\f{\log r}{r}\right)^{\f 12}$ for any smooth axially symmetric D-solutions to the Navier-Stokes equations. These improvement are based on improved weighted estimates of $\om_{\th}$ and $A_p$ weight for singular integral operators, which yields good decay estimates for $(\na u_r, \na u_z)$ and $(\om_r, \om_{z})$, where $\bm{\om}=\textit{curl }{\bf u}= \om_r {\bf e}_r + \om_{\th} {\bf e}_{\th}+ \om_z {\bf e}_z$. Another is the first decay rate estimates in the $Oz$-direction for smooth axially symmetric flows without swirl. We do not need any small assumptions on the forcing term.
\end{abstract}

\begin{center}
\begin{minipage}{5.5in}
\end{minipage}
\end{center}

\section{Introduction and main results}

In this paper, we will investigate the decay properties of the smooth axially symmetric solutions to the steady Navier-Stokes equations
\be\lab{sns} \begin{cases}
({\bf u}\cdot\nabla) {\bf u} +\nabla p -\Delta {\bf u} ={\bf f},\q \q \forall {\bf x}\in \mbR^3\\
\text{div }{\bf u}=0,\\
\displaystyle \lim_{|x|\to \oo} {\bf u}({\bf x}) ={\bf u}_{\oo}=0
\end{cases}
\ee
with finite Dirichlet integral
\be\lab{ns-dirichlet}
\int_{\mbR^3} |\nabla {\bf u}({\bf x})|^2 d{\bf x}<+\oo.
\ee
Here ${\bf u}, p$ and ${\bf f}$ denote the fluid velocity, the pressure and the body force. ${\bf u}_{\oo}$ is a constant vector and we also assume the viscosity to be $1$ for simplicity. One can also consider the same problem in the exterior domains $\Om \subset\mbR^3$ with the no-slip boundary condition on $\p\Om$, where the complement of $\Om$ is a compact axially symmetric domain. For simplicity, we only consider the whole space case and ${\bf u}_{\oo}=0$ in this paper, some of our results can also be extended to the exterior domain case.

The fundamental contribution to the existence of weak solutions to the stationary Navier-Stokes equations is due to Leray \cite{leray33}, where he constructed the weak solution of (\ref{sns}) with no-slip boundary conditions and constant velocity at infinity. Leray's solution has finite Dirichlet integral, and is usually refereed as $D$-solution. Ladyzhenskaya \cite{lady69} and Fujita \cite{fujita} also considered the nonhomogeneous boundary conditions case. It was easy to show that the $D$-solutions are smooth provided that data are smooth. In \cite{finn59,finn61}, Finn showed that any $D$-solution in three dimensional exterior domain converged uniformly pointwise to the prescribed vector ${\bf u}_{\oo}$ at infinity and, moreover, in the case ${\bf u}_{\oo}\neq 0$, he showed that if $|{\bf u}({\bf x})- {\bf u}_{\oo}|\leq C |{\bf x}|^{-\al}$ for some $\al>\f 12$ as ${\bf x}\to \oo$, then $|{\bf u}({\bf x})-{\bf u}_{\oo}|\leq \bar{C} |{\bf x}|^{-1}$ as ${\bf x}\to \oo$. Finn also suggested a calss of physical reasonable (PR) solutions to (\ref{sns}) in the three-dimensional exterior domain satisfying ${\bf u}({\bf x})= O(|{\bf x}|^{-1})$ if ${\bf u}_{\oo}=0$ or ${\bf u}({\bf x})-{\bf u}_{\oo}=O(|{\bf x}|^{-\f 12-\e})$ for some $\e>0$, if ${\bf u}_{\oo}\neq 0$. Finn \cite{finn65} then established the existence and uniqueness of a physically reasonable solution in a three dimensional exterior domain when the data are small enough. It is easy to show a PR-solution is a D-solution. However, whether the converse implication holds true has remained open for long times. In the case of ${\bf u}_{\oo}\neq 0$, Babenko \cite{babenko73} showed that every D-solution is a PR-solution if the force is of bounded support. Galdi \cite{galdi92} also proved the same result for ${\bf u}_{\oo}=0$, under the assumption that ${\bf u}$ obeys the ``energy inequality" and the viscosity is sufficiently large. In \cite{np95}, the authors established the existence and uniqueness of solution to (\ref{sns}) with the same decay rate as that of the fundamental solution of the Stokes problem, under some smallness assumptions on the data. For the investigation of the asymptotic profile of (\ref{sns}) with ${\bf u}_{\oo}=0$, one can refer to \cite{dg00,ks11,np00}. For more information about the recent results about these problems, one can refer to \cite{farwig98,fs98,galdi11}.



From the above introduction, we see that the theory about the decay properties of D-solutions to (\ref{sns}) is quite incomplete in the case of large forcing term. In this article, we will investigate the D-solution ${\bf u}$ with additional axially symmetric property to simplify this problem. On the other hand, this paper can be regarded as a continuation of my previous study \cite{cw15} with Prof. Chae, where we established some interesting Liouville type theorems for smooth axially symmetric D-solutions. From \cite{cw15}, we can see that there is a close relation between the decay properties of D-solution to (\ref{sns}) with the famous open problem of the triviality of D-solution to (\ref{sns}) with ${\bf f}\equiv 0$ and ${\bf u}_{\oo}=0$. One can also refer to \cite{chae14,cy13,cw15,galdi11,kpr15} for more recent results about this triviality problem. Now we introduce the mathematical setup of our problem. More precisely, we introduce the cylindrical coordinate
\be\no
r=\sqrt{x_1^2+x_2^2},\q \theta= \arctan\f{x_2}{x_1},\q z=x_3.
\ee
We denote ${\bf e}_r, {\bf e}_{\th}, {\bf e}_z$ the standard basis vectors in the cylindrical coordinate:
\be\no
{\bf e}_r= (\cos\th,\sin\th, 0),\q {\bf e}_{\th}= (-\sin\th,\cos\th,0),\q {\bf e}_z= (0,0,1).
\ee

A function $f$ is said to be {\it axially symmetric} if it does not depend on $\th$. A vector-valued function ${\bf u}= (u_r, u_{\th}, u_z)$ is called {\it axially symmetric} if $u_r, u_{\th}$ and $u_z$ do not depend on $\th$. A vector-valued function ${\bf u}= (u_r, u_{\th}, u_{z})$ is called {\it axially symmetric with no swirl} if $u_{\th}=0$ while $u_r$ and $u_z$ do not depend on $\th$.

Assume that ${\bf u}({\bf x})= u_r(r,z) {\bf e}_r+ u_{\th}(r,z) {\bf e}_{\th} + u_z(r,z){\bf e}_z$ is a smooth D-solution to (\ref{sns}). The corresponding asymmetric steady Navier-Stokes equations read as follows.
\be\lab{asym-sns}
\begin{cases}
(u_r\p_r+ u_z\p_z) u_r-\f{u_{\th}^2}{r} + \p_r p =\left(\p_r^2+\f{1}{r}\p_r+\p_z^2-\f{1}{r^2}\right) u_r+f_r, \\
(u_r\p_r+ u_z\p_z) u_{\th} +\f{u_r u_{\th}}{r}=\left(\p_r^2+\f{1}{r}\p_r+\p_z^2-\f{1}{r^2}\right)u_{\th}+f_{\th},\\
(u_r\p_r+ u_z\p_z) u_z + \p_z p=\left(\p_r^2+\f{1}{r}\p_r+\p_z^2\right) u_z+ f_z,\\
\p_r u_r +\f{u_r}{r} +\p_z u_z=0.
\end{cases}
\ee
Define the vorticity $\bm{\om}({\bf x})=\text{curl }{\bf u}({\bf x})= \om_r(r,z) {\bf e}_r + \om_{\th}(r,z) {\bf e}_{\th}+ \om_z(r,z) {\bf e}_z$, where
\be\no
\om_r= -\p_z u_{\th},\q \om_{\th}= \p_z u_r -\p_r u_z,\q \om_z= \f 1r \p_r(r u_{\th}).
\ee
The equations satisfied by $\om_r,\om_{\th}$ and $\om_z$ are listed as follows.
\be\lab{vorticity1}
&&(u_r\p_r+ u_z\p_z) \om_r - (\om_r\p_r +\om_z\p_z) u_r =\left(\p_r^2+\f{1}{r}\p_r+\p_z^2-\f{1}{r^2}\right)\om_r- \p_z f_{\th},\\\lab{vorticity2}
&&(u_r\p_r+ u_z\p_z) \om_{\th} - \f{u_r\om_{\th}}{r}+\f{1}{r}\p_z (u_{\th}^2) =\left(\p_r^2+\f{1}{r}\p_r+\p_z^2-\f{1}{r^2}\right)\om_{\th}+\p_z f_r- \p_r f_z,\\\lab{vorticity3}
&&(u_r\p_r+ u_z\p_z) \om_z - (\om_r\p_r +\om_z\p_z) u_z =\left(\p_r^2+\f{1}{r}\p_r+\p_z^2\right)\om_z+\f{1}{r}\p_r(r f_{\th}).
\ee

For the investigation of the decay properties of D-solutions to the exterior stationary Navier-Stokes equations (\ref{sns}), one should trace back to the important papers by Gilbarg and Weinberger \cite{gw74,gw78}. For the two dimensional exterior domain, Gilbarg and Weinberger \cite{gw74} showed that the weak solution constructed by Leray was bounded and converged to a limit ${\bf u}_0$ in a mean square sense, while the pressure converged pointwise. In \cite{gw78}, they further investigated any weak solutions to (\ref{sns}) with finite Dirichlet integral, and found that the weak solution ${\bf u}$ may not be bounded, but it grew more slowly than $(\log r)^{\f 12}$. The pressure has a finite limit at infinity and the velocity either had a limit in the mean or $\int_0^{2\pi} |{\bf u}(r,\th)|^2 d\th \to \oo$ as $r\to \oo$. The vorticity $\om= \p_{x_2} u_1-\p_{x_1} u_2$ decayed more rapidly than $r^{-\f34}(\log r)^{\f 18}$ and the first derivatives of the velocity decayed more rapidly than $r^{-\f34}(\log r)^{\f 98}$ at infinity.

Inspired by \cite{gw74} and \cite{gw78}, Choe and Jin \cite{cj09} first obtained the following decay rates for smooth axially symmetric solutions to (\ref{sns}): Let $\Om$ be an exterior domain, suppose ${\bf f}\in H^1(\Om)$ be an axially symmetric vector field with
\be\lab{cj-force}
\|{\bf f}\|_{H^{-1}(\Om)}+\|{\bf f}\|_{H^{1}(\Om)}+ \b\|\b(\f{r}{\log r}\b)^{1/2} f_{\th}\b\|_{L^2(\Om)}+\|r\p_z f_r\|_{L^2(\Om)}+ \|r\p_r f_z\|_{L^2(\Om)} \leq M
\ee
for some constant $M$. Then the axially symmetric solution $({\bf u}, p)$ with finite Dirichlet integral satisfied
\be\lab{cj1}
|u_r(r,z)|+ |u_z(r,z)|&\leq& c(M)\b(\f{\log r}{r}\b)^{\f 12},\\\lab{cj2}
|u_{\th}(r,z)| &\leq& c(M)\f{(\log r)^{1/8}}{r^{3/8}},\\\lab{cj3}
|\om_{\th}(r,z)|&\leq& \f{c(M)}{r^{\f78}}
\ee
for large $r$.

Based on these results, we can obtain better improved decay rate estimates for $u_{\th}, \bm{\om}$ and $\na {\bf u}$. We first obtain some weighted and decay estimates of $\na u_r$ and $\na u_z$ by $A_p$ weight method, since $\na u_r$ and $\na u_z$ can be expressed as singular integral operators of $\om_{\th}$. These yield weighted energy and decay estimates of $\om_r$ and $\om_z$ by using the equations of $\om_r$ and $\om_z$. Finally, we use the integral formula of $u_{\th}$ in terms of $\om_r$ and $\om_z$ to improve the decay rates. Our first main result is stated as follows.
\bt\lab{main1}{\it
Suppose ${\bf f}\in H^1(\mbR^3)$ be an axially symmetric vector field with
\be\lab{force}
\|{\bf f}\|_{H^{-1}(\mbR^3)}+\|{\bf f}\|_{H^{1}(\mbR^3)}+ \b\|\b(\f{r}{\log r}\b)^{1/2} f_{\th}\b\|_{L^2(\mbR^3)}+\|r^2 \text{curl }{\bf f}\|_{L^2(\mbR^3)}\leq M
\ee
for some constant $M$. Then the axially symmetric solution $({\bf u}, p)$ to (\ref{sns}) with finite Dirichlet integral satisfied
\be\lab{main11}
|u_{\th}(r,z)| &\leq& c(M)\b(\f{\log r}{r}\b)^{\f 12}.
\ee
for large $r$. Moreover, we have the following estimates for $\na {\bf u}$:
\be\lab{main12}
|\om_{\th}(r,z)| &\leq& c(M) r^{-(\f{19}{16})^-},\\\lab{main13}
|\na u_r(r,z)|+|\na u_z(r,z)|&\leq & c(M) r^{-(\f{9}{8})^-},\\\lab{main14}
|\om_r(r,z)|+|\om_z(r,z)|&\leq& c(M) r^{-(\f{67}{64})^-},\\\lab{main15}
|\na u_{\th}(r,z)|&\leq& c(M) r^{-(\f{67}{64})^-},
\ee
where we denote $a^-$ to be any constant less than $a$.
}
\et

Another main result in this paper is the first decay rates estimate in the $Oz$-direction for the D-solution of (\ref{sns}) without swirl. We realize that it is possible to derive the weighted energy estimates of $\Om\co \f{\om_{\th}}{r}$ in the $Oz$ direction. Together with previous weighted energy estimates on $\om_{\th}$, we can derive the decay rate of $\om_{\th}$ with respect to $\rho=\sqrt{r^2+z^2}$, which also yields the decay rate of ${\bf u}$, since ${\bf u}$ has an integral representation formula in terms of $\om_{\th}$.
\bt\lab{main2}{\it
Suppose ${\bf f}\in H^1(\Om)$ be an axially symmetric vector field without swirl, and satisfying (\ref{force}) and
\be\lab{force2}
\b\|\rho\f{(\p_z f_r-\p_r f_z)}{r}\b\|_{L^2(\mbR^3)}\leq M.
\ee

Then the axially symmetric $D$-solution $({\bf u}, p)$ to (\ref{sns}) without swirl satisfied
\be\lab{main21}
|{\bf u}(r,z)| &\leq& \f{c(M)}{(1+\rho)^{(\f{1}{8})^-}},\\\lab{main22}
|\bm{\om}(r,z)|&\leq& \f{c(M)}{(1+\rho)^{(\f{3}{8})^-}}
\ee
where $\rho=\sqrt{r^2+z^2}$.
}
\et


After this introduction section, some preliminary tools including a decay lemma and $A_p$ weight for singular integral operators will be introduced. Then we prove Theorem \ref{main1} in Section \ref{mainsection1}. We first improve the weighted energy estimates and the decay rates of $\om_{\th}$ by a bootstrap argument, these yields some good weighted estimates and decay rates of $\na u_r$ and $\na u_z$ by employing the $A_p$ weight method and the decay lemma. Based on these and the equations of $\om_r$ and $\om_z$, one can obtain some decay rates of $\om_r$ and $\om_z$, which yields better decay rates of $u_{\th}$. In Section \ref{mainsection2}, we first realizes it is possible to obtain some weighted energy estimates of $\Om$ in the $Oz$ direction, which enables us to get some decay properties of ${\bf u}$ in the $Oz$ direction.

\section{Preliminary}

\subsection{A decay lemma}

The following decay lemma is proved by the techniques developed in \cite{gw74}, \cite{gw78} and \cite{cj09}.
\bl\lab{rdecay}
{\it Suppose an smooth axially symmetric function $f(x)$ satisfies the following weighted energy estimates
\be\lab{rdecay1}
&&\int_{\mbR^3} r^{e_1} |f(r,z)|^2 dx \leq C,\\\lab{rdecay2}
&&\int_{\mbR^3} r^{e_2} |\na f(r,z)|^2 dx \leq C,\\\lab{rdecay3}
&&\int_{\mbR^3} r^{e_3} |\nabla\p_z f(r,z)|^2dx \leq C
\ee
with nonnegative constants $e_1, e_2,e_3$. Then for any $r>0$, we have
\be\lab{rdecay4}
&&\int_{-\oo}^{\oo} |f(r,z)|^2 dz \leq C r^{-\f 12(e_1+e_2)-1},\\\lab{rdecay5}
&&\int_{-\oo}^{\oo} |\p_z f(r,z)|^2 dz \leq C r^{-\f 12(e_2+e_3)-1},\\\lab{rdecay6}
&&|f(r,z)|^2\leq C r^{-\f14(e_1+2e_2+e_3)-1},\q \forall z\in \mbR^3.
\ee
}
\el

\bpf
For any integer $n\geq 0$, from (\ref{rdecay1}), we have
\be\no
\int_{2^n}^{2^{n+1}} \int_{-\oo}^{\oo} r^{e_1} |f(r,z)|^2 rdr dz \leq C.
\ee
By the intermediate value theorem, there exists $r_n\in [2^n, 2^{n+1}]$ such that
\be\no
\int_{-\oo}^{\oo} |f(r_n,z)|^2 dz \leq C r_n^{-e_1-2}.
\ee
For $\forall r>0$, choose $r_n>r$ and then
\be\no
\int_{-\oo}^{\oo} |f(r,z)|^2dz&=& \int_{-\oo}^{\oo} |f(r_n, z)|^2dz- 2\int_{r}^{r_n} \int_{-\oo}^{\oo} f(s,z)\p_s f(s,z) ds dz \co I_1 +I_2,\\\no
|I_2|&\leq& C\b(\int_{r}^{r_n} \int_{-\oo}^{\oo} \f{|f(s,z)|^2}{s^2} sds dz\b)^{\f12} \b(\int_{r}^{r_n} \int_{-\oo}^{\oo} |\p_s f(s,z)|^2 s ds dz\b)^{\f 12} \\\no
&\leq& \f{C}{r^{\f12(e_1+e_2+2)}}\b(\int_{r}^{\oo} \int_{-\oo}^{\oo} s^{e_1} |f(s,z)|^2 sds dz\b)^{\f12} \b(\int_{r}^{\oo} \int_{-\oo}^{\oo}s^{e_2}|\p_s f(s,z)|^2 s ds dz\b)^{\f 12}\\\no &\leq& O(r^{-\f 12(e_1+e_2)-1}).
\ee
Sending $n\to \oo$, $I_1\to 0$. Hence we arrive at (\ref{rdecay4}). Similarly, we also have (\ref{rdecay5}).

To prove (\ref{rdecay6}), for $|z|\leq z_1$ for some $z_1\geq 1$, we use
\be\no
|f(r,z)|^2&=&\f{1}{2z_1} \int_{-z_1}^{z_1} |f(r,t)|^2 dt+ \b(|f(r,z)|^2-\f{1}{2z_1} \int_{-z_1}^{z_1} |f(r,t)|^2 dt\b)= J_1 +J_2.
\ee
By the mean value theorem and (\ref{rdecay4})-(\ref{rdecay5}), it follows
\be\no
|J_2|&=& (|f(r,z)|^2- |f(r,z_*)|^2) \leq \int_{-z_1}^{z_1} \b|\f{\p}{\p t}|f(r,t)|^2\b| dt \\\no
&\leq& c\b(\int_{-\oo}^{\oo}|f(r,t)|^2dt\b)^{\f 12}\b(\int_{-\oo}^{\oo}|\p_t f(r,t)|^2 dt\b)^{\f 12}\\\no
&\leq& C r^{-\f14(e_1+2e_2+e_3)-1},\\\no
J_1&\leq& \f{C}{z_1 r^{\f 12(e_1+e_2)+2}}.
\ee
Then letting $z_1$ goes to $\oo$, we obtain (\ref{rdecay6}).
\epf

\subsection{$A_p$ weight and singular integral operator}
We first give the classical definition of $A_p$ weight.
\begin{definition}\label{sio1}
Let $p\in (1,\infty)$. A real valued function $w(x)$ is said to be in $A_p$ class if it satisfies
\begin{equation}\label{sio2}
\sup_{B\subset R^3}\bigg(\frac{1}{|B|}\int_Bw(x)dx\bigg)\bigg(\frac{1}{|B|}\int_B w(x)^{-\frac{p'}{p}} dx\bigg)^{\frac{p}{p'}}<\infty,
\end{equation}
where the supremum is taken over all balls $B$ in $R^3$. Here $p'$ is the H$\ddot{o}$lder conjugate of $p$, i.e. $\frac{1}{p}+\frac{1}{p'}=1$
\end{definition}

For function $w(x)\in A_p$, we can extend the Calderon-Zygmund inequality for the singular integral operator with the integral having weight function $w(x)$.

\bt\label{sio3} (\cite{stein93} page 205.)
{\it Let $p\in (1,\infty)$. Suppose $T$ is a singular integral operator of the convolution type, and $w(x)\in A^p$. Then for $f\in L^p(R^3)$,
\begin{equation}\label{sio4}
\int_{R^3}|T f(x)|^p w(x)dx\leq C \int_{R^3}|f(x)|^p w(x)dx.
\end{equation}
}
\et

\begin{lemma}\label{sio5}
For any $p\in (1,\infty)$, the function $w(x)=r^{\al}$, where $r=\sqrt{x_1^2+x_2^2}$ and $\al\in (-2, 2(p-1))$, is in $A_p$ class.
\end{lemma}

\bpf
Similar results has been proved in \cite{cl02}. Let $x_0\in \mbR^3$ be given. Set $B=B_s(x_0)$ and $d=\sqrt{x_{01}^2 + x_{02}^2}$. If $d\geq 2s$, then $d-s\leq \sqrt{x_1^2+x_2^2}\leq d+s$ for all $x\in B$. Thus if $\al\geq 0$
\begin{eqnarray}\nonumber
&\quad&\bigg(\frac{1}{|B|}\int_B w(x)dx\bigg)\bigg(\frac{1}{|B|}\int_B w(x)^{-\frac{p'}{p}} dx\bigg)^{\frac{p}{p'}}\\\nonumber
&\leq&\bigg(\frac{3}{4\pi s^3}\int_B (d+s)^{\al}dx\bigg)\bigg(\frac{3}{4\pi s^3}\int_B (d-s)^{-\frac{\al}{p-1}}dx\bigg)^{(p-1)}\\\nonumber
&=&C_2\frac{(d+s)^{\al}}{(d-s)^{\al}}\leq C_3.
\end{eqnarray}
For $\al<0$, one can also obtain a similar thing. If $d<2s$, then the cylinder $\{(x_1,x_2,x_3)\in \mbR^3|x_1^2+x_2^2<(d+s)^2, |x_3-x_3^0|<s\}$ contains the ball $B$. Thus if $\al\in (-2, 2(p-1))$,
\begin{eqnarray}\nonumber
&\quad \bigg(\frac{1}{|B|}\int_B w(x)dx\bigg)\bigg(\frac{1}{|B|}\int_B w(x)^{-\frac{p'}{p}} dx\bigg)^{\frac{p}{p'}}\\\nonumber
&\leq & \bigg(\frac{3}{4\pi s^3}2\pi \int_{x_3^0-s}^{x_3^0+s}\int_0^{d+s}\rho^{\al}\rho d\rho dx_3\bigg)\times \bigg(\frac{3}{4\pi s^3}2\pi \int_{x_3^0-s}^{x_3^0+s}\int_0^{d+s}\rho^{-\frac{\al}{p-1}+1} d\rho dx_3\bigg)^{(p-1)}\\\nonumber
&\leq & \b(\frac{3}{2 s^3}\f{2}{2+\al}s(d+s)^{2+\al}\b)\bigg(\frac{3}{2 s^3}\f{2}{2-\f{\al}{p-2}}s (d+s)^{2-\frac{\al}{p-1}}\bigg)^{(p-1)}\\\nonumber
&\leq & C(\al,p)<\oo.
\end{eqnarray}
Since the ball $B$ is arbitrary, the proof is completed.
\epf

\section{Proof of Theorem \ref{main1}.}\lab{mainsection1}

\subsection{Improved estimates for $\om_{\th}$.}

Since ${\bf f}\in H^{-1}(\mbR^3)$, the standard existence theory tells us that there exists a weak solution ${\bf u}$ to (\ref{sns}) with finite Dirichlet integral. If ${\bf f}\in H^1(\mbR^3)$, by the $L^q$ estimates of the Stokes system, then $\na^2 {\bf u}\in L^2(\mbR^3)$ and ${\bf u}\in L^{\oo}(\mbR^3)$.We first give the following basic decay estimates for $D$-solutions to (\ref{sns}).
\bl\lab{first}
{\it Let $({\bf u}, p)$ be an axially symmetric smooth D-solutions to (\ref{sns}) with the forcing term ${\bf f}$ satisfying (\ref{force}). Then there is a constant $c(M)$ such that
\be\lab{first1}
\int_{-\oo}^{\oo} |u_r(r,z)|^2 + |u_{\th}(r,z)|^2 dz &\leq& c(M),\\\lab{first2}
\int_{-\oo}^{\oo} |\om_{\th}(r,z)|^2 dz &\leq& \f{c(M)}{r}.
\ee
}
\el

\bpf
The inequality (\ref{first1}) was proven in Lemma 3.2 of \cite{cj09}. 
The second estimate (\ref{first2}) follows from (\ref{rdecay4}) in Lemma \ref{rdecay}, since $\om_{\th}$ and $\na\om_{\th}$ belong to $L^2(\mbR^3)$.
\epf

\bl\lab{omega-theta}
{\it Let $({\bf u}, p)$ be an axially symmetric smooth D-solutions to (\ref{sns}) with the forcing term ${\bf f}$ satisfying (\ref{force}). Suppose that
\be\lab{omega-theta1}
|u_r(r,z)|+ |u_{\th}(r,z)|+ |u_z(r,z)|\leq C (1+r)^{-\de},\q \forall z\in \mbR
\ee
holds for some $\de\in [0,1]$. Then the following estimates holds
\be\lab{omega-theta2}
\int_{\mbR^3}|\om_{\th}|^2 dx &\leq& c(M), \\\lab{omega-theta3}
\int_{\mbR^3}r^{1+\de}|\na\om_{\th}|^2 dx &\leq& c(M),\\\lab{omega-theta4}
\int_{\mbR^3}r^{1+3\de}|\p_z\na\om_{\th}|^2 dx &\leq& c(M).
\ee
In particular, by Lemma \ref{rdecay}, we obtain the following decay rate for $\om_{\th}$:
\be\lab{omega-theta8}
|\om_{\th}(r,z)| \leq C(M) (1+r)^{-\f 78-\f 58\de}.
\ee
}
\el

\bpf
Since $\na {\bf u}\in L^2(\mbR^3)$, (\ref{omega-theta2}) holds. Take a cut-off function $\eta\in C_0^{\oo}(\mbR^+), \eta=\eta(r)$,  satisfying $0\leq \eta\leq 1$, $\eta=0$ on $r\leq r_0$ or $r\geq 2 r_1$, $\eta=1$ on $2 r_0<r<r_1$, so that $|\na\eta|\leq \f{c}{r_0}$ on $r_0\leq r \leq 2 r_0$, $|\na \eta|\leq \f{C}{r_1}$ on $r_1<r<2 r_1$, and $\na \eta=0$ elsewhere. Here $r_0$ is a fixed positive constant. We will let $r_1$ tends to $\oo$ at the end.
Multiplying (\ref{vorticity2}) by $\eta^2 r^{a_1} \om_{\th}$ and integrating over $\mbR^3$, after some computations we obtain
\be\no
0&=& \int_{\mbR^3}\eta^2 r^{a_1}|\na \om_{\th}|^2 dx+ \int_{\mbR^3}\eta^2 r^{a_2} \b|\f{\om_{\th}}{r}\b|^2 dx -\f12\int_{\mbR^3} \p_r^2(\eta^2 r^{a_1})|\om_{\th}|^2 dx\\\no
&\q&- \f12 \int_{\mbR^3} r^{-1}\p_r(\eta^2 r^{a_1}) |\om_{\th}|^2 dx- \int_{\mbR^3} \eta^2 r^{a_2}\om_{\th}(\p_z f_r -\p_r f_z) dx -\f12\int_{\mbR^3} \p_r(\eta^2 r^{a_1}) u_r |\om_{\th}|^2 dx\\\no
&\q&- \int_{\mbR^3} \eta^2 r^{a_1-1} u_r \om_{\th}^2 dx - 2\int_{\mbR^3} \eta^2 r^{a_1-1}\om_{\th} u_{\th} \p_z u_{\th} dx \\\no
&\co& \int_{\mbR^3}\eta^2 r^{a_1}|\na \om_{\th}|^2 dx+ \int_{\mbR^3}\eta^2 r^{a_2} \b|\f{\om_{\th}}{r}\b|^2 dx + \sum_{k=1}^3 J_k + \sum_{k=1}^3 I_k.
\ee

Then for $a_1=1+\de$, we have the following estimates
\be\no
\sum_{j=1}^2|I_j|&\leq&C\int_{\mbR^3} |u_r| r^{a_1-1} \om_{\th}^2 dx \leq C\|r^{a_1-1}u_r\|_{L^{\oo}} \|\om_{\th}\|_{L^2}^2\leq C\|\om_{\th}\|_{L^2}^2,\\\no
|I_3|&\leq&C\|r^{a_1-1}u_{\th}\|_{L^{\oo}} \|\p_z u_{\th}\|_{L^2} \|\om_{\th}\|_{L^2}\leq C\|\p_z u_{\th}\|_{L^2} \|\om_{\th}\|_{L^2},\\\no
\sum_{k=1}^2 |J_k|&\leq& \int_{r\geq r_0} r^{a_1-2} |\om_{\th}|^2 dx,\q \textit{since }\q 0\leq a_1\leq 2, \\\no
|J_3|&\leq&C\|\om_{\th}\|_{L^2(\mbR^3)}\|r^2(\p_z f_r-\p_r f_z)\|_{L^2},\q\textit{since }\q 0\leq a_1\leq 2.
\ee
Assuming (\ref{omega-theta1}), letting $r_1\to \oo$, then one obtains (\ref{omega-theta3}).

Furthermore, we can obtain (\ref{omega-theta4}) by using the equation for $\p_z\om_{\th}$:
\be\lab{omega-theta-partialz}
\p_z\b((u_r\p_r+ u_z\p_z)\om_{\th}-\f{u_r}{r} \om_{\th} -\f{1}{r}\p_z({u_{\th}^2})\b)=\b(\p_r^2+\f{1}{r}\p_r+\p_z^2-\f 1{r^2}\b)\p_z\om_{\th}+\p_z(\p_z f_r-\p_r f_z).
\ee
Multiplying (\ref{omega-theta-partialz}) by $\eta^2 r^{b_1} \p_z\om_{\th}$ and integrating over $\mbR^3$, we get an integral identity with left and right hand sides as

\be\no
0&=& \int_{\mbR^3}\eta^2 r^{b_1}|\na \p_z\om_{\th}|^2 dx+ \int\eta^2 r^{b_2} \b|\f{\p_z\om_{\th}}{r}\b|^2 dx -\f12\int_{\mbR^3} \p_r^2(\eta^2 r^{b_1})|\p_z\om_{\th}|^2 dx\\\no
&\q&- \f12 \int_{\mbR^3} r^{-1}\p_r(\eta^2 r^{b_1}) |\p_z\om_{\th}|^2 dx+ \int_{\mbR^3} \eta^2 r^{b_1}\p_z\om_{\th}(\p_z f_r -\p_r f_z) dx-\int_{\mbR^3} \eta^2 r^{b_1} \p_z^2\om_{\th}[(u_r\p_r+u_z\p_z)\om_{\th}] dx \\\no
&\q&+ \int_{\mbR^3} \eta^2 r^{b_1}\p_z^2\om_{\th}\cdot\f{u_r}{r}\om_{\th} dx+ 2\int_{\mbR^3} \eta^2 r^{b_1}\p_z^2\om_{\th} \f{u_{\th}}{r}\p_z u_{\th} dx\\\no
&\co& \int_{\mbR^3}\eta^2 r^{b_1}|\na \p_z\om_{\th}|^2 dx+ \int\eta^2 r^{b_2} \b|\f{\p_z\om_{\th}}{r}\b|^2 dx+ \sum_{k=1}^3 J_k' + \sum_{k=1}^3 I_k'.
\ee
We estimate these terms as follows.
\be\no
\sum_{k=1}^2 |J_k'|&\leq&C\int_{r\geq r_0} r^{b_1-2} |\na\om_{\th}|^2 dx\leq C\int_{\mbR^3} r^{1+\de} |\na\om_{\th}|^2dx<\oo,\q \text{if }b_1-2\leq 1+\de,\\\no
|J_3'|&\leq&\int_{\mbR^3} r^{1+\de}|\na\om_{\th}|^2 dx + \int_{\mbR^3} r^{2b_1-(1+\de)}|\p_z f_r-\p_r f_z|^2 dx,\\\no
|I_1'|&\leq&\f18\int \eta^2 r^{b_1}|\na\p_z \om_{\th}|^2 dx + C \int_{\mbR^3} \eta^2 r^{b_1}|(u_r,u_z)|^2 |\na\om_{\th}|^2 dx\\\no
&\leq&\f18\int_{\mbR^3} \eta^2 r^{b_1}|\na\p_z \om_{\th}|^2 dx + C \int_{\mbR^3}  r^{1+\de}|\na\om_{\th}|^2 dx,\q \text{if }b_1\leq 1+3\de,\\\no
|I_2'|&\leq&\f18\int \eta^2 r^{b_1}|\na\p_z \om_{\th}|^2 dx + C \int_{\mbR^3} \eta^2 r^{b_1-2}|u_r|^2|\om_{\th}|^2 dx\\\no
&\leq&\f18\int_{\mbR^3} \eta^2 r^{b_1}|\na\p_z \om_{\th}|^2 dx + C \int_{\mbR^3} \eta^2 |\om_{\th}|^2 dx,\text{if }b_1-2\leq 2\de,\\\no
|I_3'|&\leq&\f18\int \eta^2 r^{b_1}|\na\p_z \om_{\th}|^2 dx + C\int_{\mbR^3} \eta^2 r^{b_1-2}|u_{\th}|^2 |\na u_{\th}|^2dx \\\no
&\leq& \f18\int_{\mbR^3} \eta^2 r^{b_1}|\na\p_z \om_{\th}|^2 dx+ C\int_{\mbR^3} |\na u_{\th}|^2 dx,\q \text{if }b_1-2\leq 2\de.
\ee

Then it is easy to see the essential restriction on $b_2$ is $b_2\leq 1+3\de$. Taking $b_2=1+3\de$ and letting $\rho_1\to \oo$, then we obtain (\ref{omega-theta4}).

\epf

\br\lab{omege-theta-remark}
{\it By (\ref{cj1})-(\ref{cj2}), we see that (\ref{omega-theta1}) holds for any $0\leq \de<\f{3}{8}$. Hence (\ref{omega-theta8}) implies that
\be\lab{omega-theta5}
|\om_{\th}(r,z)| \leq C(M) (1+r)^{-(\f{71}{64})^-}.
\ee
}
\er

\br\lab{exterior1}
{\it The case $\de=0$ was obtained in \cite{cj09}. We can also extend these arguments to the exterior domain case by choosing large enough $r_0$ in the definition of the cut-off function $\eta$.
}
\er

\bl\lab{grad-u-rz}
{\it Let $({\bf u}, p)$ be an axially symmetric smooth D-solutions to (\ref{sns}) with the forcing term ${\bf f}$ satisfying (\ref{force}). Suppose that
\be\lab{grad-u-rz1}
|u_r(r,z)|+ |u_{\th}(r,z)|+ |u_z(r,z)|\leq C (1+r)^{-\de},\q \forall z\in \mbR
\ee
holds for some $\de\in [0,1]$. Then the following estimates holds
\be\lab{grad-u-rz2}
\int_{\mbR^3}(|\na u_r|^2+|\na u_z|^2) dx &\leq& c(M), \\\lab{grad-u-rz3}
\int_{\mbR^3}r^{1+\de_1}(|\na^2 u_r|^2+ |\na^2 u_z|^2)dx &\leq& c(M),\\\lab{grad-u-rz4}
\int_{\mbR^3}r^{1+\de_2}(|\p_z\na^2u_r|^2+|\p_z \na^2 u_z|^2)dx &\leq& c(M),
\ee
where $\de_1=\de$ if $\de\in [0,1 )$ and $\de_1=1^-$ if $\de=1$, $\de_2=3\de$ if $\de<\f 13$ and $\de_2=1^-$ if $\de\in [\f 13, 1]$.

In particular, by Lemma \ref{rdecay}, we obtain the following decay rate for $\om_{\th}$:
\be\lab{grad-u-rz5}
|\na u_r (r,z)|+|\na u_z(r,z)| \leq C(M) (1+r)^{-\f 78-\f {2\de_1+\de_2}{8}}.
\ee
}\el

\bpf
Since $-\Delta (u_r {\bf e}_r + u_z {\bf e}_z)= \text{curl }(\om_{\theta} {\bf e}_{\th})$, then
\be\no
\nabla (u_r {\bf e}_r + u_z {\bf e}_z)= \nabla (-\Delta)^{-1} \text{curl }(\om_{\theta} {\bf e}_{\th}).
\ee
That is to say, $\na u_r$ and $\na u_z$ can be expressed as singular integral operators of $\om_{\th}$, and we apply Lemma \ref{sio5} to complete the proof.
%

\epf

\br\lab{grad-u-rz6}
{\it By (\ref{cj1})-(\ref{cj2}), we take $\de=(\f 38)^-$, then we have
\be\lab{grad-u-rz7}
|\na u_r (r,z)|+|\na u_z(r,z)| \leq C(M)(1+r)^{-(\f{35}{32})^-}.
\ee
}
\er

\br\lab{pressure}
{\it It seems difficult to derive weighted estimates of $(u_r, u_z)$ by using the equations of $u_r$ and $u_z$ directly, since we do not have good estimates on the pressure.
}
\er

\subsection{Estimates for $\om_r$ and $\om_z$}

\bl\lab{omega-rz}
{\it Let $({\bf u}, p)$ be an axially symmetric smooth D-solutions to (\ref{sns}) with the forcing term ${\bf f}$ satisfying (\ref{force}). Suppose that
\be\lab{omega-rz1}
|u_r(r,z)|+ |u_z(r,z)|&\leq& C (1+r)^{-\de},\\\lab{omega-rz1'}
|\na u_r(r,z)|+ |\na u_z(r,z)|&\leq& C (1+r)^{-1-\ga}
\ee
holds for some $\de,\ga\in [0,1]$. Then the following estimates holds
\be\lab{omega-rz2}
\int_{\mbR^3}(|\om_{r}|^2+|\om_z|^2) dx &\leq& c(M), \\\lab{omega-rz3}
\int_{\mbR^3}r^{1+\de\wedge\ga}(|\na\om_{r}|^2+|\na\om_z|^2) dx &\leq& c(M),\\\lab{omega-rz4}
\int_{\mbR^3}r^{1+\de\wedge\ga+2\de}(|\p_z\na\om_{r}|^2+|\p_z\na \om_z|^2) dx &\leq& c(M),
\ee
where $\de\wedge\ga= \min\{\de,\ga\}$. In particular, by Lemma \ref{rdecay}, we obtain the following decay rate for $\om_{\th}$:
\be\lab{omega-rz5}
|\om_{r}(r,z)|+|\om_z(r,z)| \leq C(M) (1+r)^{-\f 78-\f 18(3(\de\wedge\ga)+2\de)}.
\ee
}
\el
\br\lab{omega-rz6}
{\it By (\ref{cj1})-(\ref{cj2}) and (\ref{grad-u-rz7}), we take $\de=(\f{3}{8})^-$ and $\ga=(\f{3}{32})^-$ in Lemma \ref{omega-rz}, then
\be\lab{omega-rz7}
|\om_{r}(r,z)|+|\om_z(r,z)| \leq C(M) (1+r)^{-(\f{257}{256})^-}.
\ee
}\er

\bpf
We use the same cut-off function $\eta$ as in Lemma  Multiplying (\ref{vorticity1}) and (\ref{vorticity3}) by $\eta^2 r^{a_2} \om_{r}$ and $\eta^2 r^{a_2}\om_z$ respectively, integrating over $\mbR^3$ and adding them together, we obtain

\be\no
0&=& \int_{\mbR^3} \eta^2 r^{a_2}(|\na\om_r|^2+|\na\om_z|^2) dx + \int_{\mbR^3} \eta^2 r^{a_2}\b|\f{\om_r}{r}\b|^2 dx- \f12\int_{\mbR^3} \p_r^2(\eta^2 r^{a_2}) (|\om_r|^2+|\om_z|^2) dx\\\no
&\q&-\f12\int_{\mbR^3} r^{-1}\p_r(\eta^2 r^{a_2})(|\om_r|^2+|\om_z|^2) dx- \int_{\mbR^3} \eta^2 r^{a_2}\b[-\om_r\p_z f_{\th}+ \f{\om_z}{r}\p_r(r f_{\th})\b] dx\\\no
&\q&-\f12\int_{\mbR^3} \p_r(\eta^2 r^{a_2}) u_r (|\om_r|^2+|\om_z|^2) dx -\int_{\mbR^3} \eta^2 r^{a_2} [\om_r(\om_r\p_r + \om_z\p_z) u_r +\om_z (\om_r\p_r + \om_z\p_z) u_z] dx\\\no
&=& \int_{\mbR^3} \eta^2 r^{a_2}(|\na\om_r|^2+|\na\om_z|^2) dx + \int_{\mbR^3} \eta^2 r^{a_2}\b|\f{\om_r}{r}\b|^2 dx + \sum_{k=1}^3 B_k + \sum_{k=1}^2 A_k.
\ee

We can estimate $B_k, k=1,2, 3$ and $A_k, k=1,2$ as follows.
\be\no
\sum_{k=1}^2 |B_k|&\leq& \int_{r\geq r_0} r^{a_2-2} (|\om_r|^2+|\om_z|^2) dx\leq \|\om_r\|_{L^2}^2+\|\om_z\|_{L^2}^2,\q \text{if }a_2\leq 2,\\\no
|B_3|&\leq& C\|(\om_r,\om_z)\|_{L^2}(\|r^2\p_z f_{\th}\|_{L^2}+\|r\p_r (rf_{\th})\|_{L^2})<\oo,\q \text{if }a_2\leq 2,\\\no
|A_1|&\leq& C\int_{r\geq r_0} r^{a_2-1}|u_r|(|\om_r|+|\om_z|)^2dx\leq C\|(\om_r,\om_z)\|_{L^2}^2<\oo,\q \text{if } a_2\leq 1+\de,\\\no
|A_2|&\leq& C\int \eta^2 r^{a_2} |(\na u_r,\na u_z)|(|\om_r|+|\om_z|)^2 dx\leq C\|(\om_r,\om_z)\|_{L^2}^2<\oo,\q \text{if } a_2\leq 1+\ga.
\ee
Therefore, we may take $a_2=1+\de\wedge \ga$, so that
\be\no
\int_{\mbR^3} \eta^2 r^{a_2}(|\na\om_r|^2+|\na\om_z|^2) dx + \int_{\mbR^3} \eta^2 r^{a_2}\b|\f{\om_r}{r}\b|^2 dx\leq C\|(\om_r,\om_z)\|_{L^2}^2 +C (\|r^2\p_z f_{\th}\|_{L^2}+\|r\p_r (rf_{\th})\|_{L^2})<\oo.
\ee
Letting $r_1\to \oo$, we have derived \eqref{omega-rz3}.

To show (\ref{omega-rz4}), we need to use the equations for $\p_z\om_r$ and $\p_z\om_z$:
\be\lab{omega-rz-partialz1}
&&\p_z\b((u_r\p_r + u_z\p_z) \om_r - (\om_r\p_r +\om_z\p_z) u_r\b)= \b(\p_r^2+\f{1}{r}\p_r+\p_z^2-\f{1}{r^2}\b)\p_z\om_r-\p_z^2 f_{\th},\\\lab{omega-rz-partialz2}
&&\p_z\b((u_r\p_r + u_z\p_z) \om_z - (\om_r\p_r +\om_z\p_z) u_z\b)= (\p_r^2+\f{1}{r}\p_r+\p_z^2)\p_z\om_z+\p_z\b[\f{1}{r}\p_{r}(r f_{\th})\b].
\ee
Multiplying (\ref{omega-rz-partialz1}) and (\ref{omega-rz-partialz2}) by $\eta^2 r^{b_2} \p_z\om_{r}$ and $\eta^2 r^{b_2}\p_z\om_z$ respectively, integrating over $\mbR^3$ and adding them together, after some calculations we get
\be\no
0&=& \int_{\mbR^3} \eta^2 r^{b_2}(|\na\p_z\om_r|^2+|\na\p_z\om_z|^2) dx + \int_{\mbR^3} \eta^2 r^{b_2}\b|\f{\p_z\om_r}{r}\b|^2 dx- \f12\int_{\mbR^3} \p_r^2(\eta^2 r^{b_2}) (|\p_z\om_r|^2+|\p_z\om_z|^2) dx\\\no
&\q&-\f12\int_{\mbR^3} r^{-1}\p_r(\eta^2 r^{b_2})(|\p_z\om_r|^2+|\p_z\om_z|^2) dx- \int_{\mbR^3} \eta^2 r^{b_2}\b[-\p_z\om_r^2\p_z f_{\th}+ \p_z^2\om_z\b[\f{1}{r}\p_r(r f_{\th})\b] dx\\\no
&\q&-\int_{\mbR^3} \eta^2 r^{b_2}[\p_z^2\om_r (u_r\p_r+u_z\p_z)\om_r+\p_z^2\om_z (u_r\p_r+u_z\p_z)\om_z] dx \\\no
&\q&+\int_{\mbR^3} \eta^2 r^{b_2}[\p_z^2\om_r (\om_r\p_r+\om_z\p_z)u_r+\p_z^2\om_z (\om_r\p_r+\om_z\p_z)u_z] dx\\\no
&=& \int_{\mbR^3} \eta^2 r^{a_2}(|\na\p_z\om_r|^2+|\na\p_z\om_z|^2) dx + \int_{\mbR^3} \eta^2 r^{a_2}\b|\f{\p_z\om_r}{r}\b|^2 dx + \sum_{k=1}^3 B_k' + \sum_{k=1}^2 A_k'.
\ee

We estimates these terms as follows.
\be\no
\sum_{k=1}^2|B_k'|&\leq& \int_{r\geq r_0} r^{b_2-2}(|\na\om_r|+|\na\om_z|)^2 dx\leq \int_{\mbR^3} r^{1+\de\wedge\ga}(|\na \om_r|+|\na\om_z|)^2dx<\oo,\q \text{if }b_2\leq 3+\de\wedge\ga,\\\no
|B_3'|&\leq& \f18\int_{\mbR^3} \eta^2 r^{b_2} (|\na\p_z\om_r|^2+ |\na\p_z\om_z|^2) dx + C\int_{\mbR^3} r^{b_2}\b[|\p_z f_{\th}|^2+\b|\f1r\p_r(r f_{\th})\b|^2\b]^2 dx,\\\no
|A_1'|&\leq& \f18\int_{\mbR^3} \eta^2 r^{b_2} (|\na\p_z\om_r|^2+ |\na\p_z\om_z|^2) dx+ C\int_{\mbR^3} \eta^2 r^{b_2}(|u_r|+|u_z|)^2(|\na \om_r|+|\na\om_z|)^2 dx\\\no
&\leq& \f18\int_{\mbR^3} \eta^2 r^{b_2} (|\na\p_z\om_r|^2+ |\na\p_z\om_z|^2) dx+ C\int_{\mbR^3} r^{1+\de\wedge\ga}(|\na \om_r|+|\na\om_r|)^2 dx,\text{if }b_2\leq 1+\de\wedge\ga+ 2\de,\\\no
|A_2'|&\leq& \f18\int_{\mbR^3} \eta^2 r^{b_2} (|\na\p_z\om_r|^2+ |\na\p_z\om_z|^2) dx+ C\int_{\mbR^3} \eta^2 r^{b_2}(|\na u_r|+|\na u_z|)^2(|\om_r|+|\om_z|)^2 dx\\\no
&\leq& \f18\int_{\mbR^3} \eta^2 r^{b_2} (|\na\p_z\om_r|^2+ |\na\p_z\om_z|^2) dx+ C\int_{\mbR^3}(|\om_r|+|\om_r|)^2 dx,\q\text{if }b_2\leq 2+2 \ga.
\ee
Therefore, the essential restriction on $b_2$ is $b_2\leq 1+\de\wedge\ga+ 2\de$, hence we take $b_2=1+\de\wedge \ga+ 2\de$, then we infer
\be\no
\int_{\mbR^3} \eta^2 r^{b_2}(|\na\p_z\om_r|^2+|\na\p_z\om_z|^2) dx\leq \int_{\mbR^3} r^{1+\de\wedge\ga}(|\na \om_r|+|\na\om_z|)^2dx+  C\int_{\mbR^3} r^{3}\b[|\p_z f_{\th}|^2+\b|\f1r\p_r(r f_{\th})\b|^2\b]^2 dx<\oo.
\ee
Letting $r_1\to \oo$, we derive $\int_{r\geq r_0} r^{1+\de\wedge\ga+2\de} (|\na\p_z\om_r|+|\na\p_z\om_z|)^2dx<\oo$. This directly yields \eqref{omega-rz4}.

\epf

\subsection{Improved decay estimates on $u_{\th}$}

Since $\textit{curl }(u_{\th}{\bf e}_{\th})= \om_r {\bf e}_r+ \om_z {\bf e}_z$ and $\textit{div }(u_{\th}{\bf e}_{\th})=0$, then
\be\no
-\Delta( u_{\th} {\bf e}_{\th})= \textit{curl }(\om_r {\bf e}_r+ \om_z {\bf e}_z).
\ee
Fix $(r,z)\in \mbR_+\times \mbR$, we choose a smooth cut-off function $\psi\in C_0^{\oo}(\mbR^3)$, which is axially symmetric and satisfies $0\leq \psi\leq 1$,
\be\no
\psi(\rho,\ka)=\begin{cases}
1,\q r/2<\rho <2r,\q -z_1<\ka<z_1,\\
0,\q \rho<r/4, \rho>4r \ \textit{or }\  |\ka|>2 z_1,
\end{cases}
\ee
where $z_1$ is any constant such that $z_1>\max(2|z|,1)$. Since
\be\no
-\Delta(\psi u_{\th}{\bf e}_{\th})= \psi \textit{curl }(\om_r {\bf e}_r+\om_z {\bf e}_z)- 2\na\psi\cdot\na(u_{\th}{\bf e}_{\th})- (\Delta \psi) u_{\th}{\bf e}_{\th},
\ee
we get the integral representation for $u_{\th}$ in terms of $(\om_{r}, \om_z)$ and the fundamental solution $\Ga({\bf x}, {\bf y})=\Ga({\bf x}-{\bf y})=\f{1}{4\pi |{\bf x}-{\bf y}|}$ of the Laplace operator: for ${\bf x}= (r\cos\th, r\sin\th, z)$
\be\lab{u-theta-formula}\begin{array}{ll}
(\psi u_{\th}{\bf e}_{\th})(x)&= \int_{\mbR^3}\Ga(x-y)\psi(y)\textit{curl }(\om_{\rho}{\bf e}_{\rho}+ \om_{\ka}{\bf e}_{\ka}) dy-2\int_{\mbR^3}\Ga(x-y)\na\psi\cdot\na(u_{\phi}{\bf e}_{\phi})(y) dy\\
&\q-\int_{\mbR^3}\Ga(x-y)(\De\psi)(y)(u_{\phi}{\bf e}_{\phi})(y) dy\\
&=-\int_{\mbR^3}\na_y \Ga(x-y)\times[\psi(y)(\om_{\rho}{\bf e}_{\rho}+\om_{\ka}{\bf e}_{\ka})(y)] dy\\
&\q-\int_{\mbR^3}\Ga(x-y)(\na_y\psi)(y)\times (\om_{\rho}{\bf e}_{\rho}+\om_{\ka}{\bf e}_{\ka})(y)] dy\\
&\q+2\int_{\mbR^3}[\na_y\Ga(x-y)\cdot \na_y\psi(y)](u_{\phi}{\bf e}_{\phi})(y) dy+ \int_{\mbR^3}\Ga(x-y)\Delta_y\psi(y)(u_{\phi}{\bf e}_{\phi})(y)dy.
\end{array}\ee
By taking inner product to (\ref{u-theta-formula}) by ${\bf e}_{\th}$, we get the following integral representation for $u_{\th}$: for ${\bf x}= (r\cos\th, r\sin\th, z)$
\be\lab{c31}\begin{array}{ll}
u_{\th}(r,z)&=-\int_{\mbR^3} \f{\p\hat{\Ga}}{\p \ka}\psi(y)\om_{\rho}(y)\cos(\th-\phi)dy+\int_{\mbR^3} \f{\p\hat{\Ga}}{\p \rho}\psi(y)\om_{\ka}(y)\cos(\th-\phi) dy\\
&\q+ \int_{\mbR^3} \f{1}{\rho}\f{\p\hat{\Ga}}{\p \phi}\phi\om_{\ka}\sin(\th-\phi)d y-\int_{\mbR^3}  \hat{\Ga}\f{\p\psi}{\p \ka}\om_{\rho}\cos(\th-\phi) dy\\
&\q+ \int_{\mbR^3}\hat{\Ga} \f{\p\psi}{\p\rho}\om_{\ka} \cos(\th-\phi) dy+\int_{\mbR^3}\hat{\Ga} \f{1}{\rho}\f{\p\psi}{\p\phi}\om_{\ka}\sin(\th-\phi) dy\\
&\q+ 2\int_{\mbR^3}\b(\f{\p\hat{\Ga}}{\p\rho}\f{\p\psi}{\p\rho}+\f{\p\hat{\Ga}}{\p\ka}\f{\p\psi}{\p\ka}\b)u_{\phi}\cos(\th-\phi)dy+\int_{\mbR^3}\hat{\Ga}\De \psi u_{\phi} \cos(\th-\phi) dy.
\end{array}\ee

The cylindrical coordinate representation of $\Ga$ is denoted by $\hat{\Ga}= \hat{\Ga}(r,\rho, \th-\phi, z-\ka)$:
\be\no
\hat{\Ga}= \f{1}{4\pi\sqrt{r^2+\rho^2- 2r \rho \cos (\th-\phi)+ (z-\ka)^2}}.
\ee

Since ${\bf e}_{\rho}$ and ${\bf e}_r$ are different, they cause extra computations involving with $\cos(\phi-\th)$. Direct computations yield that
\be\no
\f{\p\hat{\Ga}}{\p\rho} &=&-\f{1}{4\pi} \f{\rho -r \cos(\th-\phi)}{(r^2+\rho^2- 2r\rho \cos(\th-\phi)+ (z-\ka)^2)^{\f32}},\\\no
\f{\p\hat{\Ga}}{\p \ka} &=&\f{1}{4\pi} \f{(z-\ka)}{(r^2+\rho^2- 2r\rho \cos(\th-\phi)+ (z-\ka)^2)^{\f32}}.
\ee


Define $\Ga_i= \hat{\Ga}_i(r,\rho, z-\ka)$, $i=0,\cdots, 5$, by
\be\no\begin{array}{ll}
\Ga_0= \int_{0}^{2\pi} \hat{\Ga}(r,\rho,\phi, z-\ka) d\phi,\q &\Ga_1= \int_{0}^{2\pi} \hat{\Ga}(r,\rho,\phi, z-\ka) \cos\phi d\phi,\\\no
\Ga_2= \int_{0}^{2\pi} \f{\p\hat{\Ga}}{\p \rho}(r,\rho,\phi, z-\ka) d\phi,\q &\Ga_3= \int_{0}^{2\pi} \f{\p\hat{\Ga}}{\p \rho}(r,\rho,\phi, z-\ka) \cos\phi d\phi,\\\no
\Ga_4= \int_{0}^{2\pi} \f{\p\hat{\Ga}}{\p \ka}(r,\rho,\phi, z-\ka) d\phi,\q &\Ga_5= \int_{0}^{2\pi} \f{\p\hat{\Ga}}{\p \ka}(r,\rho,\phi, z-\ka) \cos\phi d\phi.
\end{array}\ee

\bl\lab{cl22}
{\it We have the following integral representation of $u_{\th}$ in terms of $\om_r$ and $\om_z$:
\be\lab{u-theta-formula-final}\begin{array}{ll}
u_{\th}(r,z)&=-\int_{-\oo}^{\oo} \int_0^{\oo}\Ga_5\psi \om_{\rho}\rho d\rho d\ka+ \int_{-\oo}^{\oo} \int_0^{\oo}\Ga_3\psi \om_{\ka}\rho d\rho d\ka-\int_{-\oo}^{\oo} \int_0^{\oo}\Ga_1 \f{\p\psi}{\p\ka}\om_{\rho}\rho d\rho d\ka\\
&\q+ \int_{-\oo}^{\oo} \int_0^{\oo}\Ga_1\f{\p\psi}{\p \rho}\om_{\ka} \rho d\rho d\ka+ 2\int_{-\oo}^{\oo} \int_0^{\oo} (\Ga_3\f{\p\psi}{\p\rho}+\Ga_5 \f{\p\psi}{\p\ka}) u_{\phi} \rho d\rho d\ka\\
&\q+ \int_{-\oo}^{\oo} \int_0^{\oo} \Ga_1\Delta \psi u_{\phi}\rho d\rho d\ka.
\end{array}\ee
}
\el

To our purpose, we need the following estimates for $\Ga_i, i=1,\cdots, 5$.
\be\lab{c213}
|\Ga_i(r,\rho,z-\ka)|&\leq& \f{1}{\sqrt{(r-\rho)^2 + (z-\ka)^2}}\q \q \textit{for }i=0,1,\\\lab{c214}
|\Ga_i(r,\rho,z-\ka)|&\leq& \f{\rho+r}{[(r-\rho)^2 + (z-\ka)^2]^{\f32}}\q \q \textit{for }i=2,3,\\\lab{c215}
|\Ga_i(r,\rho,z-\ka)|&\leq& \f{|z-w|}{[(r-\rho)^2+(z-\ka)^2]^{\f32}}\q\q \textit{for }i=4,5.
\ee
For $\Ga_2,\Ga_3$ and $\Ga_5$, we have extra decay in $r$.
\bl\lab{cl23}
{\it Suppose $\f 14\leq \f{\rho}{r}\leq 16$, then
\be\lab{cl231}
|\Ga_2(r,\rho,z-\ka)|&\leq& \f{c}{\rho \sqrt{(\rho-r)^2+(z-\ka)^2}},\\\lab{cl233}
|\Ga_3(r,\rho,z-\ka)|&\leq& \f{c}{\rho \sqrt{(\rho-r)^2+(z-\ka)^2}},\\\lab{cl232}
|\Ga_5(r,\rho,z-\ka)|&\leq& \f{c}{r}\f{|z-\ka|}{(\rho-r)^2+(z-\ka)^2}.
\ee
}
\el
\bpf (\ref{cl231}) and (\ref{cl232}) have been proved in Lemma 2.3 in \cite{cj09}. (\ref{cl233}) follows from the following calculation
\be\no
\Ga_3(r,\rho,z-\ka)&=&-\f{1}{4\pi}\int_0^{2\pi}\f{\rho\cos\phi-r}{(r^2+\rho^2- 2r\rho \cos\phi+ (z-\ka)^2)^{\f32}}d\phi\\\no
&\q&-\f{r}{4\pi}\int_0^{2\pi}\f{\sin^2\phi d\phi}{(r^2+\rho^2- 2r\rho \cos\phi+ (z-\ka)^2)^{\f32}}\\\no
&\co&-\Ga_2(\rho,r, z-\ka)+ rJ(r,\rho,z-\ka),\\\no
|J(r,\rho,z-\ka)|&\leq&2\int_{-1}^1 \f{dt}{[r^2+\rho^2-2r\rho t+(z-\ka)^2]^{\f32}}\\\no
&=&\f2{r\rho}\b[\f{1}{\sqrt{(r-\rho)^2+(z-\ka)^2}}-\f{1}{(r+\rho)^2+(z-\ka)^2}\b]\\\no
&\leq&\f{2}{r\rho \sqrt{(r-\rho)^2+(z-\ka)^2}}.
\ee

\epf

\bl\lab{cl41}
{\it There is a positive constant $c(M)$ such that
\be\lab{c41}
|u_{\th}(r,z)|\leq c(M) \b(\f{\log r}{r}\b)^{\f12}
\ee
for large $r$, uniformly in $z$.
}\el

\bpf
From Lemma \ref{cl22}. we have
\be\no\begin{array}{ll}
u_{\th}(r,z)&=\int_{-\oo}^{\oo} \int_0^{\oo}\Ga_3\psi \om_{\ka}\rho d\rho d\ka-\int_{-\oo}^{\oo} \int_0^{\oo}\Ga_5\psi \om_{\rho}\rho d\rho d\ka+ 2\int_{-\oo}^{\oo} \int_0^{\oo} \Ga_3\f{\p\psi}{\p\rho} u_{\phi} \rho d\rho d\ka\\
&\q+2\int_{-\oo}^{\oo} \int_0^{\oo}\Ga_5 \f{\p\psi}{\p\ka}u_{\phi} \rho d\rho d\ka + \int_{-\oo}^{\oo} \int_0^{\oo}\Ga_1\f{\p\psi}{\p \rho}\om_{\ka} \rho d\rho d\ka-\int_{-\oo}^{\oo} \int_0^{\oo}\Ga_1 \f{\p\psi}{\p\ka}\om_{\rho}\rho d\rho d\ka\\\no
&\q+ \int_{-\oo}^{\oo} \int_0^{\oo} \Ga_1\b(\p_{\rho}^2\psi+\f{1}{\rho}\p_{\rho}\psi\b) u_{\phi}\rho d\rho d\ka+ \int_{-\oo}^{\oo} \int_0^{\oo} \Ga_1\p_{\ka}^2\psi u_{\phi}\rho d\rho d\ka \co \di\sum_{k=1}^8 T_k.
\end{array}
\ee

(1){\it Estimates of $T_1, T_2$.} The estimates of $T_1$ and $T_2$ are similar. We define the integration domain $D$ by
\be\no
D=\{r/4<\rho<4r, -z_1<\ka<z_1\}.
\ee
We decompose the domain of $D$ into two subregions. Let
\be\no
A&=&\{(\rho,\ka)\in [r/4,4r]\times [-z_1,z_1]: |r-\rho|\leq 1\},\\\no
B&=&\{(\rho,\ka)\in [r/4,4r]\times [-z_1,z_1]: |r-\rho|> 1\}.
\ee
Then
\be\no
T_1&\leq& r\sup_{\substack{r/4<\rho<4r\\-z_1<\ka<z_1}}  |\om_{\ka}(\rho,\ka)\int_{A} |\Ga_3(r,\rho,z-\ka)|d\rho d\ka\\\no
&\q& + \b(\int_{r/4}^{4r}\int_{-z_1}^{z_1} |\om_{\ka}(\rho,\ka)|^2\rho d\rho d\ka\b)^{\f 12} \b(\int_{B}|\Ga_3(r,\rho,z-\ka)|^2 \rho d\rho d\ka\b)^{\f 12}.
\ee

From Lemma \ref{cl23}, we have the following estimates
\be\no
&\q&\int_A |\Ga_3(r,\rho,z-\ka)| d\rho d\ka \leq \f{c}{r} \int_0^1\int_0^{z_1}\f{ds dt}{\sqrt{s^2+t^2}} \\\no
&\leq&\f{c}{r}\int_0^1 dt\b(\int_0^1 \f{d\th}{\sqrt{1+\th^2}}+ \int_1^{\f{z_1}{t}} \f{d\th}{\th}\b) \\\no
&\leq&\f{c}{r}\b(c+\int_0^1(\log z_1-\log t) dt\b)\leq \f{c+c\log z_1}{r},\\\no
&\q&\int_B |\Ga_3(r,\rho,z-\ka)|^2 d\rho d\ka \leq \f{c}{r} \int_1^{5r} \int_0^{z_1} \f{ds dt}{t^2+s^2} \\\no
&\leq& \f{c}{r} \int_1^{5r} \f{d t}{t} \int_0^{z_1/t} \f{d\th}{1+\th^2}\leq \f{c\log r}{r}.
\ee
Hence by (\ref{omega-rz7}), we have
\be\no
|T_1|
&\leq& \f{c+ c\log z_1}{r}+c\b(\f{\log r}{r}\b)^{\f 12}.
\ee

(2){\it Estimates of $T_3, T_5$ and $T_7$.} Since $|\f{\p\psi}{\p\rho}|\leq \f{c}{\rho}$,
\be\no
|T_3|
&\leq&\f{c}{r}\|u_{\phi}\|_{L^6(\mbR^3)}\b(\int_{r/4}^{r/2}\int_{-\infty}^{\infty}|\Ga_3|^{\f65}dx+\int_{2r}^{4r}\int_{-\infty}^{\infty}|\Ga_3|^{\f65}dx\b)^{\f56}\\\no
&\leq&\f{c}{r}\b[\b(\int_{r/4}^{r/2}+\int_{2r}^{4r} \b)\rho d\rho\int_{-\infty}^{\infty}\f{d\ka}{r^{6/5}[(\rho-r)^2+(z-\ka)^2]^{\f35}}\b]^{\f65}\\\no
&\leq&\f{c}{r^{1/2}}.
\ee
Note that
\be\no
\int_{-\oo}^{\oo} |\Ga_1(r,\rho,z-\ka)|^2 d\ka \leq \f{c}{|r-\rho|},
\ee
then we have
\be\no
|T_5|
&\leq&\f{c}{r}\|\om_{\ka}\|_{L^2(\mbR^3)}\b[\b(\int_{r/4}^{r/2}+\int_{2r}^{4r}\b)\rho d\rho\int_{-\infty}^{\infty}|\Ga_1|^{2}d\ka\b]^{\f12}\\\no
&\leq&\f{c}{r}\b[\b(\int_{r/4}^{r/2}+\int_{2r}^{4r} \b)\f{\rho d\rho}{|r-\rho|}\b]^{1/2}\leq \f{c}{r^{1/2}},\\\no
|T_7|
&\leq&\f{c}{r^2}\b(\int_{r/4}^{r/2}+\int_{2r}^{4r}\b)\rho d\rho \b(\int_{-\oo}^{\oo} |\Ga_1|^2 d\ka\b)^{1/2}\b(\int_{-\oo}^{\oo}|\om_{\rho}|^2 d\ka\b)^{1/2}\\\no
&\leq& \f{c}{r^2}\b(\int_{r/4}^{r/2}+\int_{2r}^{4r}\b) \f{\rho d\rho}{\sqrt{|r-\rho|}}\leq \f{c}{r^{1/2}}.
\ee
(3){\it Estimates of $T_4, T_6$ and $T_8$.} Since $|\f{\p \psi}{\p w}|\leq \f{c}{z_1}$ and $|\f{\p^2 \psi}{\p w^2}|\leq \f{c}{z_1^2}$, we observe that
\be\no
|T_4|
&\leq&\f{c}{z_1}\|u_{\phi}\|_{L^6(\mbR^3)}\b(\int_{r/4}^{4r}\int_{-\oo}^{\oo}\f{c}{r^{\f65}}\f{|z-\ka|^{\f65}}{[(r-\rho)^2+(z-\ka)^2]^{\f 35}} d\ka\b)^{\f 56}\\\no
&\leq&\f{c}{r z_1}\b(\int_{r/4}^{4r}\f{1}{|r-\rho|^{\f15}}\rho d\rho\b)^{\f 56}\leq \f{c r^{1/2}}{z_1},\\\no
|T_6|
&\leq&\f{c}{z_1}\int_{r/4}^{4r} \b(\int_{-\oo}^{\oo} |\Ga_1|^2 d\ka\b)^{1/2}\b(\int_{-\oo}^{\oo}|\om_{\rho}|^2 d\ka\b)^{1/2} \rho d\rho\\\no
&\leq&\f{c}{z_1}\int_{r/4}^{4r}\f{1}{\sqrt{|r-\rho|}} \f{1}{\sqrt{\rho}}\rho d\rho\leq \f{c r}{z_1},\\\no
|T_8|
&\leq&\f{c}{z_1^2}\int_{r/4}^{4r} \b(\int_{-\oo}^{\oo} |\Ga_1|^2 d\ka\b)^{1/2}\b(\int_{-\oo}^{\oo}|u_{\phi}|^2 d\ka\b)^{1/2} \rho d\rho\\\no
&\leq&\f{c}{z_1^2}\int_{r/4}^{4r}\f{1}{\sqrt{|r-\rho|}} \rho d\rho\leq \f{c r^{3/2}}{z_1^2}.
\ee

We choose $z_1=r^2$, then combining all the above estimates, we finally obtain (\ref{c41}).

%
%
%

\epf

\bpf[Proof of Theorem \ref{main1}.]
We have proved (\ref{main11}) in Lemma \ref{cl41}. By (\ref{cj1}) and (\ref{main11}), we see that (\ref{omega-theta1}) holds for $\de=(\f 12)^-$, hence by (\ref{omega-theta8}), we have (\ref{main12}). Applying Lemma \ref{grad-u-rz}, we can take $\de_1=(\f 12)^-$ and $\de_2=1^-$, hence by (\ref{grad-u-rz5}), we obtain (\ref{main13}). Then we use Lemma \ref{omega-rz}, where we can take $\de=(\f 12)^-, \ga=(\f 18)^-$, (\ref{omega-rz5}) implies (\ref{main14}). Since $\na u_{\th}$ can be expressed as singular integral operators of $\om_r$ and $\om_{z}$, we can use $A_p$ weight to derive same weighted energy estimates of $\na u_{\th}$ from (\ref{omega-rz2})-(\ref{omega-rz4}). Then by Lemma \ref{rdecay}, we obtain (\ref{main15}).

\epf

\section{Proof of Theorem \ref{main2}.}\lab{mainsection2}
\bpf[Proof of Theorem \ref{main2}.]
{\bf Step 1.} We have the following weighted estimates for $\Om\co \f{\om_{\th}}{r}$.

{\bf Claim 1.} {\it Suppose that
\be\lab{urz100}
|u_r(r,z)|+ |u_z(r,z)|\leq C (1+\rho)^{-\tau},\q \rho=\sqrt{r^2+z^2}
\ee
for some $\tau\in [0,1]$, then we have
\be\lab{Omega1}
&&\int_{\mbR^3} \rho^{1+\tau} |\nabla\Om(r,z)|^2 rdr dz <\oo,\\\lab{Omega2}
&&\int_{\mbR^3} \rho^{1+3\tau} |\nabla\p_z\Om(r,z)|^2 rdr dz <\oo.
\ee
}
To prove {\bf Claim 1.}, we see that for the smooth axially symmetric flows with no swirl, $\Om$ satisfies the following equation
\be\lab{Omega-no swirl}
(u_r\p_r+ u_z\p_z)\Om = (\p_r^2+\f{3}{r}\p_r+ \p_z^2)\Om+\f{1}{r}(\p_z f_r- \p_r f_z).
\ee
Choose a cut-off function $\phi\in C_0^{\oo}(\mbR^3)$, $\phi=\phi(\rho)$, satisfying $0\leq \phi\leq 1$, $\phi(\rho)=1$ on $2\rho_0\leq\rho \leq \rho_1$, $\phi(\rho)=0$ on $\rho\leq \rho_0$ or $\rho\geq 2\rho_1$, such that $|\na \phi|\leq \f{C}{\rho_0}$ on $\rho_0<\rho<2\rho_0$ and $|\na \phi|\leq \f{C}{\rho_1}$ on $\rho_1<\rho< 2\rho_1$, and $\na \phi=0$ elsewhere.

Multiplying (\ref{Omega-no swirl}) by $\phi^2 \rho^{d_1} \Om$ and integrating over $\mbR^3$, then we get
\be\no
0&=& \int \phi^2 \rho^{d_1} |\na\Om|^2 dx -\f12\int \p_r(\phi^2 \rho^{d_1}) |\Om|^2 dx-\f12 \int \p_z^2(\phi^2 \rho^{d_1}) |\Om|^2 dx \\\no
&=& 3\int \phi\phi' \rho^{d_1-1} |\Om|^2 dx -\int \phi^2 \rho^{d_1}\Om \f1r(\p_z f_r-\p_r f_z) dx-\f12 \int_{\mbR^3} \b[u_r\p_r(\phi^2 \rho^{d_1})+ u_z\p_z(\phi^2 \rho^{d_1})\b] |\Om|^2 dx\\\no
&=&\int \phi^2 \rho^{d_1} |\na\Om|^2 dx + \sum_{i=1}^3 K_{1i} + L_{11}.
\ee
We estimate these terms as follows.
\be\no
\sum_{i=1}^3 |K_{1i}|&\leq& C\int_{\rho\geq \rho_0}  \rho^{d_1-2} |\Om|^2 dx\leq \|\Om\|_{L^2}^2,\q\text{if }d_1\leq 2,\\\no
|L_{11}|&\leq& \int \phi |(u_r,u_z)|\rho^{d_1-1} |\Om|^2 dx\leq \|\Om\|_{L^2}^2,\q\text{if }d_1\leq 1+\tau.
\ee
This yields \eqref{Omega1} by letting $\rho_1\to \oo$.

%

To derive the estimate (\ref{Omega2}), we use the equation for $\p_z\Om$.
\be\lab{partialz-Omega}
\p_z\b((u_r\p_r+u_z\p_z) \Om\b) =(\p_r^2+\f{3}{r}\p_r + \p_z^2) \p_z\Om+\f{1}{r}\p_z(\p_z f_r- \p_r f_z).
\ee
Multiplying (\ref{partialz-Omega}) by $\phi^2 \rho^{d_2} \p_z\Om$ and integrating over $\mbR^3$, then we get
\be\no
0&=& \int \phi^2 \rho^{d_1}|\na\p_z\Om|^2 dx -\f12 \int [\p_r^2(\phi^2 \rho^{d_2})+\p_z^2(\phi^2 \rho^{d_2})] |\p_z\Om|^2 dx+\f{3}{2}\int r^{-1}\p_r(\phi^2) \rho^{d_2} |\p_z\Om|^2 dx\\\no
&\q&+ \int \phi^2 \rho^{d_2}\b(\p_r^2\Om+\f{\p_r\Om}{r}\b)\b[\f1r(\p_z f_r-\p_r f_z)\b] dx + \int \p_r(\phi^2 \rho^{d_2})\p_r\Om \b[\f1r(\p_z f_r-\p_r f_z)\b] dx\\\no
&\q&-\int \phi^2 \rho^{d_2}\b(\p_r^2\Om+\f{\p_r\Om}{r}\b)(u_r\p_r+u_z\p_z)\Om dx - \int \p_r(\phi^2 \rho^{d_2})\p_r\Om (u_r\p_r+ u_z\p_z)\Om dx\\\no
&=&\int \phi^2 \rho^{d_1}|\na\p_z\Om|^2 dx + \sum_{i=1}^5 K_{2i} + \sum_{j=1}^2 L_{2j}.
\ee
These terms can be bounded as follows.
\be\no
\sum_{i=1}^3 |K_{2i}|&\leq& \int_{\rho\geq \rho_0} \rho^{d_2-2} |\na \Om|^2 dx<\oo,\q\text{if } d_2\leq 3+\tau,\\\no
|K_{24}|&\leq&\f18 \int \phi^2 \rho^{d_2}|\na \p_z\Om|^2 dx + C\int \phi^2 \rho^{d_2}\b[\f1r(\p_z f_r-\p_r f_z)\b]^2 dx,\\\no
|K_{25}|&\leq& \int_{\rho\geq \rho_0} \rho^{d_2-2} |\na\Om|^2 dx + C\int \phi^2 \rho^{d_2}\b[\f1r(\p_z f_r-\p_r f_z)\b]^2 dx,\q \text{if }d_2\leq 3+\tau,\\\no
|L_{21}|&\leq& \f18 \int \phi^2 \rho^{d_2}|\na \p_z\Om|^2 dx + C\int \phi^2 \rho^{d_2}|(u_r,u_z)|^2 |\na\Om|^2 dx,\q \text{if } d_2\leq 1+3\tau,\\\no
|L_{22}|&\leq& \int_{\rho\geq \rho_0} \rho^{d_2-1} |(u_r, u_z)||\na \Om|^2 dx<\oo,\q \text{if }d_2\leq 2(1+\tau).
\ee
From the above estimates, we derive \eqref{Omega2} by letting $\rho_1\to \oo$.

{\it Step 2. Now we can derive the decay rate for $\om_{\th}$.}

{\bf Claim 2.} {\it Suppose that \eqref{urz100} holds, we can infer that
\be\lab{omega1}
|\om(r,z)|\leq C(1+\rho)^{-(\f{5}{16}+\f{1}{2}\tau)^-}.
\ee}

Now we prove {\bf Claim 2.} Combining the results in Lemma \ref{omega-theta} and (\ref{Omega1})-(\ref{Omega2}), then
\be\lab{Omega3}
&&\int_{\mbR^3} r^2 |\Om(r,z)|^2 dx<\oo,\\\lab{Omega4}
&&\int_{\mbR^3} (r^{3+\de}+|z|^{1+\tau}) |\na\Om(r,z)|^2 dx<\oo,\\\lab{Omega5}
&&\int_{\mbR^3} (r^{3+3\de}+|z|^{1+3\tau}) |\na\p_z\Om(r,z)|^2 dx<\oo,
\ee
where $\de$ can be any constant less than $\f 12$. Fix $d>1$, then for each $n\in \mb{N}$,
\be\no
\int_{2^n}^{2^{n+1}} \int_d^{\oo} r^2 |\Om(r,z)|^2 rdr dz<\oo.
\ee
By mean value theorem, there exists $z_n\in [2^n, 2^{n+1}]$ such that
\be\no
\int_d^{\oo} r^2 |\Om(r,z_n)|^2 rdr \leq \f{C}{z_n}.
\ee
Then for any $z$, choose $z_n>z$ and
\be\no
\int_d^{\oo} |\Om(r,z)|^2 r dr &=&\int_d^{\oo} |\Om(r,z_n)|^2 r dr- 2\int_d^{\oo} \int_z^{z_n}\Om(r,t)\p_t \Om(r,t) r dr dt\co I_1 +I_2,\\\no
|I_2|&\leq&\b(\int_d^{\oo}\int_z^{z_n} |\Om(r,t)|^2 r drdt\b)^{1/2}\b(\int_d^{\oo}\int_z^{z_n} |\p_t\Om(r,t)|^2 r drdt\b)^{1/2}\leq \f{C}{d |z|^{\f12(1+\tau)}}.
\ee
Letting $z_n\to \oo$, then $I_1\to 0$ and
\be\lab{Omega6}
\int_d^{\oo}|\Om(r,z)|^2 rdr\leq \f{C}{d|z|^{\f12(1+\tau)}}.
\ee
Similarly, one can find $z_n\in [2^n, 2^{n+1}]$ such that
\be\no
\int_d^{\oo} |\na\Om(r,z_n)|^2 rdr &\leq& \f{C}{z_n^2},\\\no
\int_d^{\oo} |\na\Om(r,z)|^2 r dr&=&\int_d^{\oo} |\na\Om(r,z_n)|^2 r dr- 2\int_d^{\oo}\int_z^{z_n} \na\Om(r,t)\cdot \p_t\na\Om(r,t) r dr dt\co J_1+J_2,\\\no
|J_2|&\leq&\b(\int_d^{\oo}\int_z^{z_n} |\na\Om(r,t)|^2 r drdt\b)^{1/2}\b(\int_d^{\oo}\int_z^{z_n} |\p_t\na\Om(r,t)|^2 r drdt\b)^{1/2}\\\no
&\leq&\begin{cases}
\f{C}{d^{3+2\de}},\\
\f{C}{|z|^{1+2\tau}}.
\end{cases}
\ee
Letting $n\to \oo$, $J_1\to 0$. Take $\de=(\f 12)^-$ and $J_2\leq \min\{\f{C}{d^{3+2\de}}, \f{C}{|z|^{1+2\tau}}\}$
\be\no
\int_d^{\oo} |\na\Om(r,z)|^2 r dr\leq \b(\f{C}{d^{4^-}}\b)^{\f14} \b(\f{C}{|z|^{1+2\tau}}\b)^{\f{3}{4}}\leq \f{C}{d |z|^{(\f 34(1+2\tau))^-}}.
\ee

Finally,
\be\no
|\Om(d,z)|^2 &=&\f{1}{r_1-d}\int_d^{r_1} |\Om(r,z)|^2 dr+ (|\Om(r,z)|^2-\f{1}{r_1-d}\int_d^{r_1}|\Om(r,z)|^2 dr)\co H_1 +H_2,\\\no
|H_2|&=&\b||\Omega(d,z)|^2-|\Om(d_*,z)|^2\b| \leq 2\int_d^{r_1} |\Om(r,z)\p_r\Om(r,z)|dr\\\no
&\leq&\f{C}{d}\b(\int_d^{\oo}|\Om(r,z)|^2 rdr\b)^{1/2}\b(\int_d^{\oo}|\na \Om(r,z)|^2 r dr\b)^{1/2}\\\no
&\leq&\f{C}{d}\b(\f{C}{d|z|^{\f12(1+\tau)}}\b)^{1/2}\b(\f{C}{d |z|^{(\f 34(1+2\tau))^-}}\b)^{\f 12}\leq \f{C}{d^2|z|^{(\f58+\tau)^-}},
\ee
which implies that
\be\lab{omega-z}
|\om(d,z)|\leq \f{C}{|z|^{(\f5{16}+\f12\tau)^-}}.
\ee
Together with Theorem \ref{main1}, we have
\be\lab{omega-decay}
|\om(r,z)|\leq \f{C}{\rho^{(\f5{16}+\f12\tau)^-}},\q \forall (r,z)\in \mbR_+\times \mbR, \rho=\sqrt{r^2+z^2}.
\ee

{\bf Step 3.} Now we derive the new decay rate of ${\bf u}$. Fix any ${\bf x}\in\mbR^3\setminus\{0\}$, define a cut-off function $\psi\in C_0^{\oo}(\mbR^3)$ satisfying $\psi({\bf y})\equiv 1$ for $\forall {\bf y}\in B_{\rho/4}({\bf x})$ and $\psi({\bf y})\equiv 0$ for $\forall {\bf y}\not\in B_{\rho/2}({\bf x})$, where $\rho=|{\bf x}|$. One can require that $|\na\psi({\bf y})|\leq \f{C}{|y|}, |\na^2 \psi({\bf y})|\leq \f{C}{|y|^2}$ for $\forall {\bf y}\in D\co B_{\rho/2}({\bf x})\setminus B_{\rho/4}({\bf x})$. Setting ${\bf u}({\bf x})= u_r {\bf e}_r +u_z {\bf e}_z$, since $\text{curl }{\bf v}= \om_{\th} {\bf e}_{\th}$, then
\be\lab{v-omega}\begin{array}{ll}
{\bf v}({\bf x}) &= -\int_{\mbR^3} \na_{{\bf y}}\Ga({\bf x},{\bf y}) \times (\om_{\phi}({\bf y})\psi({\bf y}){\bf e}_{\phi}) d {\bf y} -\int_{\mbR^3} \Ga({\bf x},{\bf y})((\na_{{\bf y}}\psi({\bf y})\times {\bf e}_{\phi}))\om_{\phi}({\bf y}) d {\bf y}\\
&\quad+ \int_{\mbR^3} \Ga({\bf x},{\bf y})(\Delta_{{\bf y}}\psi)({\bf y}) {\bf v}({\bf y}) d{\bf y} + 2 \int_{\mbR^3} (\na_{{\bf y}}\Ga)({\bf x},{\bf y})\cdot (\na_{{\bf y}}\psi)({\bf y}) {\bf v}({\bf y}) d {\bf y}\\
&\co K_1 +K_2 +K_3 +K_4.
\end{array}\ee
We estimate $K_i, i=2,3,4$ as follows.
\be\no
|K_2|&\leq& \f{C}{\rho}\b(\int_D |\Ga({\bf x}-{\bf y})|^2 d{\bf y}\b)^{1/2}\b(\int_D |\om_{\phi}({\bf y})|^2 d {\bf y}\b)^{1/2}\leq \f{C}{\rho^{1/2}},\\\no
|K_3|&\leq& \f{C}{\rho^2}\b(\int_D |\Ga({\bf x}-{\bf y})|^{\f 65} d{\bf y}\b)^{\f 56}\b(\int_D |{\bf v}({\bf y})|^6 d {\bf y}\b)^{\f 16}\leq \f{C}{\rho^{1/2}},\\\no
|K_4|&\leq& \f{C}{\rho}\b(\int_D |\na \Ga({\bf x}-{\bf y})|^{\f 65} d{\bf y}\b)^{\f 56}\b(\int_D |{\bf v}({\bf y})|^6 d {\bf y}\b)^{\f 16}\leq \f{C}{\rho^{1/2}}.
\ee
For the estimate of $K_1$, fix a $d\in (0,\f{\rho}{2})$, which will be determined later, then
\be\no
|K_1|&\leq& \sup_{{\bf y}\in B_d({\bf x})}|\om_{\phi}({\bf y})|\int_{B_d({\bf x})}|\na\Ga({\bf x}-{\bf y})| d{\bf y}\\\no
&\q&+ \b(\int_{B_{\rho/2}({\bf x})\setminus B_{d}({\bf x})} |\na\Ga({\bf x}-{\bf y})|^2 d{\bf y}\b)^{\f 12} \b(\int_{B_{\rho/2}({\bf x})\setminus B_{d}({\bf x})}|\om_{\phi}({\bf y})|^2 d {\bf y}\b)^{\f 12}\\\no
&\leq& C\rho^{-(\f{5}{16}+\f12\tau)^-} d+ C d^{-\f 12}.
\ee
By choosing $d=\rho^{(\f{5}{24}+\f13\tau)^-}$, we obtain the optimal bound for $|K_1|\leq \f{C}{\rho^{(\f5{48}+\f16\tau)^-}}$. Hence we have
\be\lab{urz200}
|(u_r,u_z)(r,z)|\leq C (\rho+1)^{-(\f5{48}+\f16 \tau)^-}
\ee

{\bf Step 4. Iteration.} At the beginning, we have $\tau=0$ in \eqref{urz100}, then by using the arguments developed in {\bf Step 1} to {\bf Step 3}, we have a new $\tau$ in \eqref{urz100}, which will be denoted by $\tau_1=(\f{5}{48})^-$. Run a second iteration of these three steps, we get a new $\tau_2= \tau_1+ \f16 \tau_1$, and after $n$ iteration, we get
\be\no
\tau_n= \tau_{1}+ \f16 \tau_{n-1}= \tau_1 \sum_{i=0}^{n-1}\f1{6^i}.
\ee
Let $n\to \oo$, $\tau_n\to (\f18)^-$ as $n\to \oo$. In a word, we infer the following decay rates
\be\lab{urz300}
|(u_r,u_z)(r,z)|&\leq& C (\rho+1)^{-(\f18)^-},\\\lab{omgea300}
|\om(r,z)| &\leq& C(\rho+1)^{-(\f38)^-}.
\ee

\epf

\br\lab{gamma}
{\it Since there are no improved decay estimates of $\Ga_2,\Ga_3,\Ga_5$ in $w$ in (\ref{cl231})-(\ref{cl233}), it seems difficult to use the argument in Lemma \ref{cl41} to get better decay estimates of $u_{\th}$. So we use (\ref{v-omega}) directly. From the estimates of $K_2, K_3$ and $K_4$, one may also see the difficulties to improve the decay rates in Theorem \ref{main1}.
}
\er

\br\lab{exterior2}
{\it
One can also extend Theorem \ref{main2} to the exterior domain case.
}
\er

{\bf Acknowledgement.}  The author would like to thank Prof. Dongho Chae and Prof. Zhouping Xin for the stimulating discussions and constant encouragement and supports. Special thanks also go to the referee for the important suggestions and comments, which make this paper more readable.



\end{document}